\documentclass[12pt]{article}
\usepackage{array}
\usepackage{amsthm,amsmath,amssymb,color}
\usepackage{fullpage}
\usepackage{graphicx}
\usepackage{float}
\usepackage{tikz}
\usetikzlibrary{shapes.geometric,calc}

\newtheorem{thm}{Theorem}[section]

\newtheorem{lem}[thm]{Lemma}

\newtheorem{cor}[thm]{Corollary}

\theoremstyle{remark}

\theoremstyle{plain}

\newtheorem{claim}{Claim}

\newcommand{\Gbar}{\overline{G}}
\DeclareMathOperator{\rank}{rank}
\DeclareMathOperator{\col}{col}

\begin{document}

\title{Co-edge-regular graphs with two $\lambda$-values\\ and four distinct eigenvalues}
\author{Suil O\thanks{Department of Applied Mathematics and Statistics, The State University of New York, Korea, Incheon, 21985, suil.o@sunykorea.ac.kr. Research supported by the National Research Foundation of Korea (NRF) grant funded by the Korea government(MSIT) No. RS-2025-23523950.}
}

\maketitle

\begin{abstract}
\noindent
A connected $k$-regular graph $G$ on $n$ vertices is co-edge-regular if any two distinct non-adjacent vertices have the same number of common neighbors. For adjacent vertices $u$ and $v$, let $\lambda(u,v)=|N(u)\cap N(v)|$. Ge and Koolen asked whether the theorem of Neumaier on strongly regular graphs with fixed smallest eigenvalue extends to co-edge-regular graphs with exactly four distinct eigenvalues in which $\lambda(u,v)$ takes exactly two values over the edges $uv$. In this paper, we answer this question. For every integer $t\ge2$ there is $n_t$ such that if $G$ is such a graph with smallest eigenvalue at least $-t$, $n\ge n_t$, and $k<n-2-(t-1)^2/4$, then either $G$ is a clique extension of a strongly regular graph, or $G$ is obtained from a strongly regular graph with a Hoffman coloring by adding all edges inside the color classes. Consequently, $G$ is a clique extension of a strongly regular graph, a $p\times q$ grid with $p\ne q$, or the block graph of a resolvable Steiner $2$-design in which every parallel class is made into a clique.
\\


\end{abstract}

\section{Introduction}

For undefined terms in graph theory, see West~\cite{W}; for the basics of spectral graph theory, see Brouwer and Haemers~\cite{BH} or Godsil and Royle~\cite{GR}.

Let $G$ be a connected non-complete $k$-regular graph on $n$ vertices. For adjacent vertices $u$ and $v$, let $\lambda(u,v)=|N(u)\cap N(v)|$. The graph $G$ is \emph{co-edge-regular} with parameters $(n,k,\mu)$ if any two distinct non-adjacent vertices have exactly $\mu$ common neighbors. We say that $G$ has \emph{two $\lambda$-values} if $\lambda(u,v)$ takes exactly two values over the edges $uv$ of $G$ (in the terminology of Ge and Koolen~\cite{GK}, such a graph has level~$2$). A co-edge-regular graph whose $\lambda(u,v)$ takes a single value is exactly a strongly regular graph, and a connected regular graph has at most three distinct eigenvalues if and only if it is complete or strongly regular.

A central result on strongly regular graphs is a theorem of Sims~\cite[Theorem~8.6.4]{BvM}, proved by Neumaier~\cite{N}: for fixed $m\ge2$, all but finitely many connected strongly regular graphs with smallest eigenvalue $-m$ are Steiner graphs, Latin square graphs, or complete multipartite graphs. Tan, Koolen and Xia~\cite{TKX} proposed an analogue for four eigenvalues. They conjectured that, for a fixed lower bound on the smallest eigenvalue and under suitable largeness and density assumptions, every connected co-edge-regular graph with exactly four distinct eigenvalues is either the $b$-clique extension of a strongly regular graph or a $p\times q$ grid. Both of these families have two $\lambda$-values.

Ge and Koolen~\cite{GK} disproved this conjecture. One of their families consists of graphs with two $\lambda$-values: the block graph of a resolvable Steiner $2$-design in which every parallel class is made into a clique. Their other families, and those of van Dam, Ge and Koolen~\cite{vDGK}, have at least three $\lambda$-values; in fact, the number of $\lambda$-values is unbounded among co-edge-regular graphs with exactly four distinct eigenvalues~\cite{vDGK}, so no analogue of the theorem of Neumaier can be expected without a restriction on it, as observed by Yang, Gebremichel, Ur Rehman, Yang and Koolen~\cite{YGRYK}. Ge and Koolen asked whether the theorem of Neumaier can be generalized to co-edge-regular graphs with two $\lambda$-values and exactly four distinct eigenvalues~\cite[Problem~34(3)]{GK}.

To state the answer, we recall two constructions. For a graph $Q$ and an integer $b\ge2$, the \emph{$b$-clique extension} of $Q$ is obtained by replacing every vertex of $Q$ by a clique of order $b$ and completely joining two such cliques whenever the corresponding vertices of $Q$ are adjacent. A \emph{Hoffman coloring} of a regular graph is a partition of its vertex set into cocliques attaining the Hoffman bound (see Section~\ref{sec:def}). For a strongly regular graph $S$, where we also allow a disjoint union of at least two complete graphs of the same order, and a Hoffman coloring $\mathcal P$ of $S$, let $S^{\mathcal P}$ be the graph obtained from $S$ by joining every two vertices in the same part of $\mathcal P$. The grids and the family of Ge and Koolen are of this form.

In this paper, we answer the question of Ge and Koolen. The hypotheses of Theorem~\ref{thm:main} are exactly those of the conjecture of Tan, Koolen and Xia.

\begin{thm}\label{thm:main}
For every integer $t\ge2$ there is $n_t$ with the following property. Let $G$ be a connected co-edge-regular graph on $n\ge n_t$ vertices with valency $k<n-2-(t-1)^2/4$, two $\lambda$-values, exactly four distinct eigenvalues, and smallest eigenvalue at least $-t$. Then one of the following holds.
\begin{enumerate}
\item[(I)] $G$ is the $b$-clique extension of a connected non-complete strongly regular graph for some $b\ge2$.
\item[(II)] $G=S^{\mathcal P}$, where $S$ is a strongly regular graph, possibly a disjoint union of complete graphs of the same order, and $\mathcal P$ is a Hoffman coloring of $S$.
\end{enumerate}
\end{thm}

Combined with the theorem of Neumaier, Theorem~\ref{thm:main} gives exactly three families.

\begin{cor}\label{cor:three}
Under the hypotheses of Theorem~\ref{thm:main}, $G$ is the $b$-clique extension of a connected non-complete strongly regular graph, a $p\times q$ grid with $p\ne q$, or the block graph of a resolvable Steiner $2$-design in which every parallel class is made into a clique.
\end{cor}

The proof uses the fact that such a graph is the union of two relations of a symmetric $3$-class association scheme whose third relation is the complement. In the imprimitive case, one of the two $\lambda$-classes is a disjoint union of cliques, and the two types appear according to how these cliques are joined. In the primitive case, we show that both $\lambda$-classes carry a partial linear space with large lines, and then derive a contradiction from the spectrum of the scheme.

In Section~\ref{sec:def}, we give definitions and tools. In Section~\ref{sec:main}, we prove Theorem~\ref{thm:main} and Corollary~\ref{cor:three}. In Section~\ref{sec:remark}, we make concluding remarks.

\section{Definitions and Tools}\label{sec:def}

Throughout, $G$ is a connected non-complete $k$-regular graph on $n$ vertices with adjacency matrix $A$, and $\Gbar$ is its complement. When $G$ has exactly four distinct eigenvalues, we write them as $k=\theta_0>\theta_1>\theta_2>\theta_3$. We use $J$ for the all-ones matrix, $\mathbf 1$ for the all-ones vector, and $\chi_X$ for the characteristic vector of a set $X$.

The $b$-clique extension of a graph $Q$ has adjacency matrix $(A_Q+I)\otimes J_b-I$. If $Q$ has eigenvalue $\theta$, then its $b$-clique extension has eigenvalue $b(\theta+1)-1$, and the remaining eigenvalue is $-1$.

Let $H$ be a $d$-regular graph on $N$ vertices with smallest eigenvalue $\tau$. Every coclique of $H$ has at most $N(-\tau)/(d-\tau)$ vertices~\cite{BH}. A coclique attaining this bound is a \emph{Hoffman coclique}, and a partition of $V(H)$ into Hoffman cocliques is a \emph{Hoffman coloring}. For the graphs $S^{\mathcal P}$ of Section~1, we note the following. The $p\times q$ grid is $S^{\mathcal P}$ with $S=qK_p$ and the parts of $\mathcal P$ transversal to the components of $S$. If $S$ is the block graph of a resolvable Steiner $2$-$(m,K,1)$ design and $\mathcal P$ is a resolution, then $S^{\mathcal P}$ is the family of Ge and Koolen~\cite[Section~4]{GK}. If $S$ is connected and not complete multipartite, then the parts of $\mathcal P$ are Delsarte cliques of the primitive strongly regular graph $\overline S$, and by the theorem of Haemers and Tonchev on spreads~\cite{HT} the relations ``adjacent in $S$'', ``in the same part of $\mathcal P$'' and ``neither'' form a $3$-class association scheme. If $S=qK_p$, then these three relations are ``same row'', ``same column'' and ``neither'' on the $p\times q$ grid, which form the product of the trivial schemes on $p$ and on $q$ points. In both cases $S^{\mathcal P}$ is co-edge-regular.

A \emph{symmetric $3$-class association scheme} on a set $V$ consists of relations $R_0,R_1,R_2,R_3$ partitioning $V\times V$, with $R_0$ the diagonal, whose adjacency matrices $I,A_1,A_2,A_3$ span a commutative algebra closed under ordinary multiplication; we write $A_iA_j=\sum_hp^h_{ij}A_h$. Here $p^h_{ij}$ is the number of $w$ with $(u,w)\in R_i$ and $(w,v)\in R_j$ for $(u,v)\in R_h$, and $p^h_{ij}=p^h_{ji}$. Let $k_i=p^0_{ii}$ be the valency of $R_i$. The algebra has primitive idempotents $E_0=J/n,E_1,E_2,E_3$ with $A_j=\sum_iP_{ij}E_i$; let $m_i=\rank E_i$. Then
\begin{equation}\label{eq:mult}
m_i=n\Big/\sum_{j=0}^3P_{ij}^2/k_j
\end{equation}
by the orthogonality relations~\cite{BH}. The scheme is \emph{imprimitive} if $R_0\cup\bigcup_{i\in I}R_i$ is an equivalence relation for some $\emptyset\ne I\subsetneq\{1,2,3\}$, and \emph{primitive} otherwise; in a primitive scheme every $R_i$ is a connected graph.

A \emph{partial linear space} consists of points and lines, each line a set of points, with any two points on at most one line. Its \emph{point graph} joins two points if they lie on a common line. If $N$ is its point-line incidence matrix and every point lies on exactly $r$ lines, then $NN^{\top}=A+rI$ for the adjacency matrix $A$ of the point graph.

We use the following known results.

\begin{thm}[Neumaier~\cite{N}]\label{thm:N}
For each integer $m\ge2$, all but finitely many connected strongly regular graphs with smallest eigenvalue $-m$ are Steiner graphs, that is, block graphs of Steiner $2$-designs with blocks of size $m$, Latin square graphs $LS_m(q)$, or complete multipartite graphs.
\end{thm}

\begin{thm}[Metsch~{\cite[Result 2.2]{M}}]\label{thm:M}
Let $H$ be a $d$-regular graph, and let $m,c,e$ be positive integers. Suppose that any two adjacent vertices have exactly $c$ common neighbors, any two non-adjacent vertices have at most $e+1$ common neighbors, $c>(2m-1)e-1$, and $d<(m+1)(c+1)-\frac{m(m+1)}2e$. Call a maximal clique with at least $c+2-(m-1)e$ vertices a \emph{line}. Then every vertex lies on at most $m$ lines, and any two adjacent vertices lie on a unique common line.
\end{thm}

Metsch states this for graphs in which any two adjacent vertices have at least $\lambda^{(1)}$ and at most $\lambda^{(2)}$ common neighbors and any two non-adjacent vertices have at most $\mu$ common neighbors, with $2\lambda^{(1)}-\lambda^{(2)}>(2m-1)(\mu-1)-1$ and every vertex of degree less than $(m+1)(\lambda^{(1)}+1)-\frac{m(m+1)}2(\mu-1)$; Theorem~\ref{thm:M} is the case $\lambda^{(1)}=\lambda^{(2)}=c$ and $\mu=e+1$.

Koolen, Gebremichel, Yang and Yang~\cite{KGYY} proved that for each $t\ge2$ there is $C(t)$ such that every co-edge-regular graph with smallest eigenvalue at least $-t$ and $k\ge C(t)$ satisfies $\mu\le t^2(t-1)$ or $n-k-1\le(t-1)^2/4+1$. Finally, Delsarte, Goethals and Seidel~\cite{DGS} proved that a set of unit vectors in $\mathbb R^d$ with at most three distinct inner products between distinct members has at most $\binom{d+2}3+\binom{d+1}2$ elements.

\section{Proof of Theorem~\ref{thm:main}}\label{sec:main}

We first fix notation and record two lemmas. Let $G$ be as in Theorem~\ref{thm:main}, let $\lambda_1\ne\lambda_2$ be its two $\lambda$-values, let $R_1$ and $R_2$ be the sets of ordered pairs of adjacent vertices with $\lambda(u,v)=\lambda_1$ and $\lambda(u,v)=\lambda_2$, and let $R_3$ be the set of ordered pairs of distinct non-adjacent vertices.

\begin{lem}\label{lem:scheme}
The relations $R_0,R_1,R_2,R_3$ form a symmetric $3$-class association scheme. Its eigenspaces are the eigenspaces of $A$.
\end{lem}

\begin{proof}
We have $A=A_1+A_2$ and $A^2=kI+\lambda_1A_1+\lambda_2A_2+\mu A_3$ with $A_3=J-I-A$. Since $G$ is connected and regular, $J\in\mathbb R[A]$. Since $\lambda_1\ne\lambda_2$, these equations give $A_1,A_2,A_3\in\mathbb R[A]$. As $A$ has four distinct eigenvalues, $\dim\mathbb R[A]=4$, and the non-zero $0/1$ matrices $I,A_1,A_2,A_3$ have disjoint supports, so they form a basis of $\mathbb R[A]$. Hence $\mathbb R[A]$ is a commutative algebra spanned by $I,A_1,A_2,A_3$ and closed under multiplication. Its primitive idempotents are the projections onto the eigenspaces of $A$.
\end{proof}

From now on, $p^h_{ij}$, $k_i$, $P_{ij}$, $E_i$ and $m_i$ refer to this scheme, with $E_i$ the eigenspace of $A$ for $\theta_i$. Since $\mu\ge1$, $G$ has diameter $2$, so $n\le1+k^2$ and $k\to\infty$ as $n\to\infty$. By the result of Koolen, Gebremichel, Yang and Yang and the assumption $k<n-2-(t-1)^2/4$,
\begin{equation}\label{eq:mu}
\mu\le M_0:=t^2(t-1)\quad\text{for }n\ge n_t.
\end{equation}
Below, ``for $n\ge n_t$'' means for $n$ larger than a constant depending only on $t$, and all constants depend only on $t$.

\begin{lem}\label{lem:local}
Let $E_0'=\binom{t^2}2(\mu-1)$ and $C_0=(t^2-1)(\mu-1)+E_0'$. For every vertex $u$ there are disjoint sets $D_1,\dots,D_\alpha$ and $D_\ast$ with $\alpha\le t^2$ and $N(u)=D_1\cup\dots\cup D_\alpha\cup D_\ast$ such that every $D_i$ is a clique, $|D_\ast|\le E_0'$, every vertex of $D_i$ has at most $\mu-1$ neighbors in $D_j$ for $j\ne i$, and $|D_i|-1\le\lambda(u,v)\le|D_i|-1+C_0$ for $v\in D_i$.
\end{lem}

\begin{proof}
An induced $K_{1,j}$ has smallest eigenvalue $-\sqrt j$, so by interlacing~\cite{BH} the subgraph induced by $N(u)$ has independence number $\alpha\le t^2$. Let $\{w_1,\dots,w_\alpha\}$ be a maximum coclique of it. Let $D_\ast$ be the set of vertices of $N(u)$, other than the $w_i$, that are adjacent to at least two of the $w_i$, and let $D_i$ consist of $w_i$ and the vertices of $N(u)$ whose only neighbor among the $w_j$ is $w_i$. By maximality these sets partition $N(u)$. If two vertices of $D_i\setminus\{w_i\}$ were non-adjacent, replacing $w_i$ by them would give a larger coclique; hence $D_i$ is a clique. A vertex of $D_\ast$ is a common neighbor, other than $u$, of two non-adjacent vertices $w_i,w_j$, so $|D_\ast|\le\binom\alpha2(\mu-1)$. A vertex $v\in D_i$ is not adjacent to $w_j$ for $j\ne i$, and its neighbors in $D_j\setminus\{w_j\}$ are common neighbors of $v$ and $w_j$ other than $u$; so there are at most $\mu-1$ of them. Counting $N(u)\cap N(v)$ inside $D_i$, inside the other $D_j$, and inside $D_\ast$ gives the bounds on $\lambda(u,v)$.
\end{proof}

\begin{proof}[Proof of Theorem~\ref{thm:main}]
We treat the imprimitive and the primitive case separately.

\medskip\noindent\emph{The imprimitive case.}

Assume first that the scheme is imprimitive.

\begin{claim}\label{cl:types}
After exchanging $R_1$ and $R_2$ if necessary, $R_0\cup R_1$ is an equivalence relation whose classes, called \emph{fibres}, are cliques of $G$ of a common order $b\ge2$. Moreover, there is an integer $a\ge1$ such that every vertex has $0$ or $a$ $R_2$-neighbors in each fibre not containing it. If $a=b$, then (I) holds; if $a<b$, then every vertex has exactly $a$ $R_2$-neighbors in every fibre other than its own.
\end{claim}

\begin{proof}
If $R_0\cup R_3$ is an equivalence relation, then $\Gbar$ is a disjoint union of cliques, so $G$ is complete multipartite and has at most three distinct eigenvalues. If $R_0\cup R_1\cup R_2$ is an equivalence relation, then $G$ is disconnected. Suppose that $R_0\cup R_1\cup R_3$ is an equivalence relation, with classes $C_1,\dots,C_r$ of a common order $c$; here $r\ge2$ since $R_2\ne\emptyset$. Pairs in different classes lie in $R_2$, so every vertex has $c-1-(n-k-1)$ neighbors in its own class. For $x=\chi_{C_1}-\chi_{C_2}$ a direct count gives $Ax=-(n-k)x$. Hence $-t\le-(n-k)$, so $n-k-1\le t-1\le(t-1)^2/4+1$, contradicting $k<n-2-(t-1)^2/4$. The union $R_0\cup R_2\cup R_3$ is symmetric to this. Hence $R_0\cup R_1$ is an equivalence relation; its classes are cliques of $G$ of a common order $b\ge2$. Let $\ell\ge2$ be the number of fibres, so $n=\ell b$.

If $(u,v)\in R_2$ and $v\in Y$, then the other $R_2$-neighbors of $u$ in $Y$ are the $w$ with $(v,w)\in R_1$ and $(w,u)\in R_2$; so $a=1+p^2_{12}$.

Let $a=b$. If some $u\in X$ has an $R_2$-neighbor in $Y$, then $u$ is joined to all of $Y$, and every vertex of $Y$ is joined to $b$ vertices of $X$, hence to all of $X$. So between two fibres either all pairs or no pairs are in $R_2$. Let $Q$ be the graph on the fibres recording the first alternative. Then $A_2=A_Q\otimes J_b$ and $A_3=(J-I-A_Q)\otimes J_b$, so $A_2^2\in\operatorname{span}\{I,A_1,A_2,A_3\}$ gives $A_Q^2\in\operatorname{span}\{I,A_Q,J\}$. Since $G$ is connected and non-complete, $Q$ is connected, non-complete and strongly regular, and $G$ is its $b$-clique extension.

Let $a<b$. Take an $R_2$-pair $(u,v)$ with $u\in X$, $v\in Y$, and $u'\in X$ with $(u',v)\in R_3$. For $(x,y)\in R_3$, the number $p^3_{12}$ equals the number of $R_2$-neighbors of $y$ in the fibre of $x$; at $(u',v)$ it is $a$. If a vertex $y$ had no $R_2$-neighbor in a fibre $X'$, then any $x\in X'$ would give $(x,y)\in R_3$ with this number equal to $0$. Hence every vertex has exactly $a$ $R_2$-neighbors in every other fibre.
\end{proof}

From now on, let $a<b$. Then $R_2$ is $k_2$-regular with $k_2=a(\ell-1)$. Let $W$ be the space of vectors whose sum over every fibre is $0$. Since $A_2\chi_Y=a(\mathbf 1-\chi_Y)$ for every fibre $Y$, on the span of the vectors $\chi_X-\chi_Y$ we have $A_2=-a$ and $A=b-1-a\ge0$. On $W$ we have $A_1=-I$, $J=0$ and $A_3=-A_2$.

\begin{claim}\label{cl:W}
On $W$, $A_2$ has exactly two eigenvalues $\rho>0>-\tau$. With $\nu=a(b-a)/(b-1)$,
\begin{enumerate}
\item[(i)] $\theta_3=-\tau-1$, so $\tau\le t-1$;
\item[(ii)] $\rho\tau=(\ell-1)\nu$, $\rho-\tau=p^2_{22}-p^3_{22}$, and the multiplicity of $\rho$ is $\ell(b-1)\tau/(\rho+\tau)$;
\item[(iii)] for any two fibres $X\ne Y$, all singular values of the restriction to $\chi_X^\perp$ of the $R_2$-biadjacency matrix $B_{XY}$ are at most $\tau$; in particular $1\le\nu\le\tau^2$.
\end{enumerate}
\end{claim}

\begin{proof}
On $W$, the equation $A_2^2=k_2I+p^1_{22}A_1+p^2_{22}A_2+p^3_{22}A_3$ becomes $A_2^2=(k_2-p^1_{22})I+(p^2_{22}-p^3_{22})A_2$. If $A_2$ had only one eigenvalue on $W$, then $G$ would have at most three distinct eigenvalues. Counting the common $R_2$-neighbors of $u$ and the other vertices of its fibre gives $(b-1)p^1_{22}=k_2(a-1)$, so the product of the two eigenvalues is $-(k_2-p^1_{22})=-(\ell-1)\nu<0$. Since $b-1-a\ge0$, the smallest eigenvalue of $G$ is $-\tau-1$, which gives (i). From $\operatorname{tr}A_2=0$ we get $k_2-a(\ell-1)+m\rho-m'\tau=0$ for the multiplicities $m,m'$ of $\rho,-\tau$ on $W$; together with $m+m'=\ell(b-1)$ this gives (ii).

For (iii), $B=B_{XY}$ is $a$-regular on both sides, so its singular values are $a$ and those of its restriction to $\chi_X^\perp$. Let $x\perp\chi_X$ and $y\perp\chi_Y$ be unit vectors supported on $X$ and $Y$ with $x^{\top}By$ equal to the largest such singular value $\sigma_2$, and let $z=x-y\in W$. Since fibres are $R_2$-cocliques, $z^{\top}A_2z=-2\sigma_2\ge-\tau\|z\|^2=-2\tau$. Finally $ab=\operatorname{tr}BB^{\top}\le a^2+(b-1)\sigma_2^2$ gives $\nu\le\tau^2$, and $a(b-a)\ge b-1$ gives $\nu\ge1$.
\end{proof}

If $\tau=a$, then (II) holds. Indeed, the eigenvalues of $A_2$ are among $k_2,\rho,-a$. If $R_2$ is connected, it is strongly regular. If not, $k_2$ has multiplicity at least $2$, so $\rho=k_2$ and $R_2$ is a disjoint union of complete graphs of the same order. In both cases the smallest eigenvalue of $R_2$ is $-a$, and every fibre is an $R_2$-coclique of order $b=n\cdot a/(k_2+a)$, which is the Hoffman bound. So the fibres form a Hoffman coloring $\mathcal P$ of $S=R_2$, and $G=S^{\mathcal P}$.

\begin{claim}\label{cl:taua}
For $n\ge n_t$, we have $\tau=a$.
\end{claim}

\begin{proof}
Suppose that $\tau\ne a$, and write $P=p^3_{22}$.

\emph{Step 1.} $n\le\binom{d+2}3+\binom{d+1}2$ with $d=(b-1)(t-1)^2$, and $a\le t^2(t-1)/2$.

Let $F$ be the primitive idempotent for the eigenspace of $A_2$ on $W$ for $\rho$; this is the eigenspace of $A$ for $\rho-1$, since $\rho-1=b-1-a$ would leave $A$ with only three distinct eigenvalues. By Claim~\ref{cl:W}(ii) and $\nu\le\tau^2$, its rank satisfies $\ell(b-1)\tau^2/((\ell-1)\nu+\tau^2)\le(b-1)\tau^2/\nu\le d$. Since $F$ is in the Bose--Mesner algebra, $F=\sum_i c_iA_i$ with $c_0>0$, and the vectors $Fe_u$ have equal length and inner products $c_i$ for $(u,v)\in R_i$. Suppose $Fe_u=Fe_v$ for some $u\ne v$. Then $\{i:c_i=c_0\}$ defines a non-trivial equivalence relation. It cannot contain $R_1$, since then every vector of the eigenspace would be constant on fibres while lying in $W$. It cannot be $R_0\cup R_2$: then $R_2$ would be a disjoint union of complete graphs, so its smallest eigenvalue would be $-1$ and $\tau=1$; moreover the $a$ $R_2$-neighbors of a vertex in another fibre would be pairwise $R_2$-adjacent while lying in a common fibre, so $a=1=\tau$. It cannot be $R_0\cup R_2\cup R_3$, since two vertices of a fibre are both related by $R_2\cup R_3$ to a vertex of another fibre. Unions with $R_3$ but not $R_2$ were excluded in Claim~\ref{cl:types}. Hence the vectors $Fe_u$ are $n$ distinct vectors of equal length with at most three inner products, and the theorem of Delsarte, Goethals and Seidel gives the first bound.

For the second, let $(u,v)\in R_3$ with $u\in X$ and $v\in Y$. The common neighbors of $u$ and $v$ are the $a$ $R_2$-neighbors of $v$ in $X$, the $a$ $R_2$-neighbors of $u$ in $Y$, and $P$ further vertices. So $\mu=2a+P$, and \eqref{eq:mu} applies.

By Step~1, we may assume $b\ge B_0$ for a constant $B_0$ chosen large enough for each step below; in particular $b\ge4a$, so $a/\nu=(b-1)/(b-a)\le2$.

\emph{Step 2.} $(\ell-1)(\tau a-\nu)=\tau(a-\tau)+\tau Pb/a$.

Fix $u$ and count pairs $(w,v)$ with $w\in N_2(u)$ and $v\in N_2(w)\setminus\{u\}$, where $N_2$ denotes the set of $R_2$-neighbors, according to the relation of $v$ to $u$. This gives $k_2(k_2-1)=(b-1)p^1_{22}+k_2p^2_{22}+(b-a)(\ell-1)P$. With $(b-1)p^1_{22}=k_2(a-1)$ we get $ap^2_{22}=a(k_2-a)-(b-a)P$. By Claim~\ref{cl:W}(ii), $(p^2_{22}-P+\tau)\tau=(\ell-1)\nu$. Eliminating $p^2_{22}$ and using $k_2=a(\ell-1)$ gives the claim.

\emph{Step 3.} $P\ge1$, $\rho\ge c_1b$ for a constant $c_1>0$, and $p^1_{22}$ is bounded.

Suppose $P=0$. Then $p^2_{22}=a(\ell-2)$. For an $R_2$-pair $(u,v)$ with $v\in Y$, the common $R_2$-neighbors lie outside the fibres of $u$ and $v$, and $u$ has exactly $a(\ell-2)$ $R_2$-neighbors there; so $N_2(u)\setminus Y\subseteq N_2(v)$. If $a=1$, this gives $\{u\}\cup N_2(u)=\{v\}\cup N_2(v)$, so $R_2$ is a union of complete graphs and $\tau=1=a$. If $a\ge2$, take distinct $w,w'\in N_2(u)\cap Y$; then $p^1_{22}\ge1+a(\ell-2)$, and since $(b-1)p^1_{22}=a(a-1)(\ell-1)$ and $\ell-1\le1+a(\ell-2)$, we get $b\le a(a-1)+1<B_0$. Hence $P\ge1$.

The right side of Step~2 is then at least $-\tau^2+\tau b/a>0$, so $\tau a-\nu>0$ and $\ell-1\ge b/a^2-\tau$. Hence $\rho=(\ell-1)\nu/\tau\ge(\ell-1)/\tau\ge c_1b$. If $a=1$, then $p^1_{22}=0$. If $a\ge2$, then $\nu\ge3a/4\ge3/2$, so $\tau\ge\sqrt{3/2}$ and $\tau a-\nu\ge a(\tau-1)\ge a/5$; then Step~2 gives $\ell-1\le C_1b$, so $p^1_{22}=a(a-1)(\ell-1)/(b-1)$ is bounded.

\emph{Step 4.} $R_2$ is the point graph of a partial linear space in which every point lies on exactly $\tau$ lines; in particular $\tau$ is an integer. Moreover $a\le\tau$.

Apply Theorem~\ref{thm:M} to each connected component of $R_2$, with $d=k_2$, $c=p^2_{22}$, $e=\max\{1,\max(p^1_{22},P)-1\}$ and $m=2\lceil\tau\rceil$. The non-adjacent pairs of $R_2$ are the pairs in $R_1\cup R_3$, so they have at most $e+1$ common $R_2$-neighbors. By Claim~\ref{cl:W}(ii), $c=\rho-\tau+P\ge\rho-\tau$, which is large by Step~3, while $e$ is bounded. Also $k_2=(a/\nu)\tau\rho\le2\tau\rho$, so $(m+1)(c+1)-k_2\ge\rho-(2\tau+1)(\tau-1)$, which exceeds $\frac{m(m+1)}2e$. Hence every $R_2$-edge lies on a unique line. With $C=(m-1)e$, every line has between $c+2-C$ and $c+2$ points, the upper bound because two points of a line have the rest of it as common $R_2$-neighbors.

Let $r_u$ be the number of lines through $u$. Every $R_2$-neighbor of $u$ lies on exactly one line through $u$, so $\sum_{L\ni u}(|L|-1)=k_2$ and $r_u\in[k_2/(c+1),k_2/(c+1-C)]$. This interval does not depend on $u$ and has length at most $2\tau\rho C/((\rho-\tau)(\rho-\tau-C))<1$, so $r_u=r$ is constant. With $N$ the point-line incidence matrix, $NN^{\top}=A_2+rI$, so $-\tau\ge-r$. If $r>\tau$, then the kernel of $NN^{\top}$ is contained in the span of the vectors $\chi_X-\chi_Y$, so $\rank N\ge n-\ell+1$; but the number of lines is at most $nr/(c+2-C)\le nm/(\rho-\tau-C)<n-\ell+1$. Hence $r=\tau$. Finally, a line is an $R_2$-clique and so meets every fibre at most once; hence $|L|\le\ell$ and $a(\ell-1)=k_2\le\tau(\ell-1)$.

\emph{Step 5.} By Step~4 and $\tau\ne a$, we have $\tau>a$. Then $NN^{\top}=A_2+\tau I$ has the eigenvalue $\tau-a>0$ on the span of the vectors $\chi_X-\chi_Y$, the eigenvalues $k_2+\tau$ and $\rho+\tau$, and its kernel is the eigenspace of $-\tau$. Hence $\rank N=\ell+\ell(b-1)\tau/(\rho+\tau)$ by Claim~\ref{cl:W}(ii). On the other hand, with $C'=\tau+C$, the number of lines is at most $n\tau/(c+2-C)\le\ell b\tau/(\rho-C')$. So
\[
\frac{b\tau}{\rho-C'}\ge1+\frac{(b-1)\tau}{\rho+\tau},\quad\text{that is,}\quad\tau\bigl(\rho+b\tau+(b-1)C'\bigr)\ge(\rho-C')(\rho+\tau).
\]
For $\rho\ge4C'$ the right side is at least $\rho^2/2+\tau\rho$, while the left side is at most $\tau\rho+\tau b(\tau+C')$. Since $\rho\ge c_1b$ by Step~3, this is impossible for $b\ge B_0$.
\end{proof}

By Claims~\ref{cl:types} and~\ref{cl:taua}, (I) or (II) holds in the imprimitive case.

\medskip\noindent\emph{The primitive case.}

Assume now that the scheme is primitive, so $R_1$, $R_2$ and $R_3$ are connected. We show that this is impossible for $n\ge n_t$. Let $\delta=|\lambda_1-\lambda_2|$.

We first note that
\begin{equation}\label{eq:large}
k_1,k_2\ge(n-1)^{1/3}/2\quad\text{and}\quad k_j\le E_0'+t^2(\lambda_j+1)\ \ (j=1,2);
\end{equation}
in particular $k_1,k_2,\lambda_1,\lambda_2\to\infty$ as $n\to\infty$. Indeed, the number of $R_1$-walks of length $i$ between $u$ and $v$ is $(A_1^i)_{uv}$, which depends only on the relation of $(u,v)$. So the $R_1$-distance between two vertices depends only on their relation, and the connected graph $R_1$ has diameter at most $3$. Hence $n\le1+k_1+k_1^2+k_1^3$; similarly for $R_2$. An $R_j$-neighbor $v\notin D_\ast$ of $u$ lies in some $D_i$ of Lemma~\ref{lem:local} with $|D_i|\le\lambda_j+1$, which gives the second inequality.

We say that a partial linear space whose point graph is $G$ is \emph{monochromatic} if every line is contained in $R_1$ or in $R_2$, that is, all pairs of points on a line lie in the same relation.

\begin{claim}\label{cl:geom}
For $n\ge n_t$, $G$ is the point graph of a monochromatic partial linear space in which every point lies on exactly $r_j$ lines contained in $R_j$ for $j=1,2$. Hence, with $N_j$ the incidence matrix of the $R_j$-lines and $N=[N_1\ N_2]$, we have $A_j=N_jN_j^{\top}-r_jI$, $A=NN^{\top}-(r_1+r_2)I$, $\theta_3=-(r_1+r_2)$ and $E_3=\ker N^{\top}$.
\end{claim}

The last statements follow from the first: every edge lies on exactly one line, the number of lines is less than $n$ since lines are large, so $-(r_1+r_2)$ is the smallest eigenvalue of $A$ with eigenspace $\ker N^{\top}=\ker N_1^{\top}\cap\ker N_2^{\top}$, and by Lemma~\ref{lem:scheme} this is $E_3$. We prove the first statement separately for $\delta>C_0$ and $\delta\le C_0$.

\begin{proof}
\emph{Case $\delta>C_0$.}
Two vertices $v,v'$ of the same $D_i$ in Lemma~\ref{lem:local} have $\lambda(u,v)$ and $\lambda(u,v')$ within $C_0<\delta$ of each other, so they have the same relation to $u$. By \eqref{eq:large}, $k_1>E_0'$, so $u$ has an $R_1$-neighbor $v\notin D_\ast$; let $v\in D_i$. A common neighbor $w$ of $u$ and $v$ with $(u,w)\in R_2$ lies outside $D_i$, so there are at most $C_0$ of them. Hence $p^1_{12}+p^1_{22}\le C_0$, using $p^1_{21}=p^1_{12}$, and similarly $p^2_{12}+p^2_{11}\le C_0$.

Apply Theorem~\ref{thm:M} to each connected component of $R_1$ with $d=k_1$, $c=p^1_{11}=\lambda_1-2p^1_{12}-p^1_{22}\ge\lambda_1-2C_0$, $e=\max(C_0,M_0)$ and $m=t^2$. Non-adjacent pairs of $R_1$ lie in $R_2\cup R_3$ and have at most $\max(p^2_{11},p^3_{11})\le e+1$ common $R_1$-neighbors. By \eqref{eq:large},
\[
(m+1)(c+1)-\tfrac{m(m+1)}2e-k_1\ge\lambda_1-(t^2+1)(2C_0-1)-\tfrac{t^2(t^2+1)}2e-E_0'-t^2>0
\]
and $c>(2m-1)e-1$ for $n\ge n_t$. So every $R_1$-edge lies on a unique line, which is an $R_1$-clique with between $c+2-C$ and $c+2$ points, where $C=(m-1)e$. The number of lines through $u$ lies in $[k_1/(c+1),k_1/(c+1-C)]$, an interval independent of $u$ of length at most
\[
\frac{(E_0'+t^2(\lambda_1+1))C}{(\lambda_1-2C_0+1-C)(\lambda_1-2C_0+1)}<1,
\]
so it is a constant $r_1$. The same argument applies to $R_2$. A line of $R_1$ and a line of $R_2$ cannot share two points, so together they form a monochromatic partial linear space with point graph $G$.

\emph{Case $\delta\le C_0$.}
Let $\lambda_{\min}=\min(\lambda_1,\lambda_2)$ and $\lambda_{\max}=\lambda_{\min}+\delta$. Fix $u$ and take $D_1,\dots,D_\alpha,D_\ast$ from Lemma~\ref{lem:local}. If $v\in N(u)\setminus D_i$ is non-adjacent to some $w\in D_i$, then $N(v)\cap D_i$ and $u$ are common neighbors of $v$ and $w$, so $|N(v)\cap D_i|\le\mu-1$; otherwise $v$ is adjacent to all of $D_i$. A vertex $z\in D_\ast$ is adjacent to all of some $D_i$, since otherwise $\lambda(u,z)\le t^2(\mu-1)+E_0'<\lambda_{\min}$, and to all of at most one, since otherwise $\lambda(u,z)\ge2(\lambda_{\min}+1-C_0)>\lambda_{\max}$. Add $z$ to that $D_i$, and call the results $D_i'$. Two vertices added to the same $D_i$ have more than $\mu$ common neighbors, so each $D_i'$ is a clique; every vertex of $D_i'$ has at most $M_1=\mu-1+E_0'$ neighbors in $D_j'$ for $j\ne i$; and $\lambda_{\min}+1-C_0\le|D_i'|\le\lambda_{\max}+1$.

Let $L_i(u)=\{u\}\cup D_i'$ and $v\in D_i'$. The clique $L_i(u)\setminus\{v\}\subseteq N(v)$ has more than $t^2M_1$ vertices, so it lies in a single part $D_j'(v)$ of the decomposition at $v$: some part meets it in more than $M_1$ vertices, and a vertex of the clique in another part would be adjacent to all of them, and $L_i(u)\subseteq L_j(v)$; by symmetry $L_i(u)=L_j(v)$. So the sets $L_i(u)$ are the lines of a partial linear space with point graph $G$ in which every edge lies on exactly one line, and every line $L$ has $h_L$ points with $h_{\min}\le h_L\le h_{\max}$, where $h_{\max}-h_{\min}\le2C_0$ and $h_{\min}\ge\lambda_{\min}+2-C_0$. The number $t_u$ of lines through $u$ lies in $[k/(h_{\max}-1),k/(h_{\min}-1)]$; since $k\le t^2(h_{\max}-1)$, this interval has length less than $1$, so $t_u=r$ is constant. Hence $A=NN^{\top}-rI$ for the incidence matrix $N$, $\theta_3=-r$, and $E_3=\ker N^{\top}$.

It remains to show that every line is monochromatic. Let $\Pi=I-E_3$ be the projection onto $\col N$. Since $\Pi$ is in the Bose--Mesner algebra, it has a constant diagonal $d/n$, with $d=n-m_3=\rank N$, and constant entries $\pi_1,\pi_2,\pi_3$ on $R_1,R_2,R_3$. Write $P_{3j}=-\rho_j$, so $\rho_1+\rho_2=r$, and $\varepsilon=m_3/n$. Since $\Pi=\sum_{i\le2}E_i$ and $\sum_iP_{ij}m_i=0$ for $j\ne0$, we have $\pi_j=\varepsilon\rho_j/k_j$.

We have $\pi_1\ne\pi_2$. Otherwise $\Pi=xI+yA+zJ$ for some reals $x,y,z$; comparing the eigenvalues on $E_1$ and $E_2$ gives $y(\theta_1-\theta_2)=0$, so $y=0$ and $x=1$, and then $\Pi$ acts as $1$ on $E_3$, contradicting $\Pi E_3=0$. Since $\Pi\chi_L=\chi_L$ for every line $L$, reading this at a point $u\in L$ gives $d/n+a_1\pi_1+a_2\pi_2=1$, where $a_j$ is the number of points of $L$ in relation $R_j$ with $u$; as $a_1+a_2=h_L-1$ and $\pi_1\ne\pi_2$, the numbers $a_j=a_j(L)$ depend only on $L$.

Let $S=\rho_1^2/k_1+\rho_2^2/k_2+(r-1)^2/k_3$, so $n/m_3=1+S$ by \eqref{eq:mult}. Since $d\le\#\text{lines}\le nr/h_{\min}$, we get $S\le r/(h_{\min}-r)$ and $m_3>n/2$. So $P_{3j}$, an eigenvalue of an integer matrix of multiplicity greater than $n/2$, is rational and hence an integer. If $\rho_1<0$, then $\rho_2\ge r+1$ and $S\ge(r+1)^2/k\ge(r+1)^2/(r(h_{\max}-1))$, which exceeds $r/(h_{\min}-r)$ for $n\ge n_t$. If $\rho_1=0$, then $S\ge r^2/(k-k_1)\ge r/(h_{\max}-1-k_1/r)$, so $k_1\le r(h_{\max}-h_{\min}+r-1)$, contradicting \eqref{eq:large}. The same arguments with $1$ and $2$ exchanged exclude $\rho_2\le0$. Hence $\rho_1,\rho_2\ge1$. Now
\[
\frac{\rho_1^2}{k_1}+\frac{\rho_2^2}{k_2}=\frac{r^2}k+\frac{(\rho_1k_2-\rho_2k_1)^2}{k_1k_2k}
\quad\text{and}\quad
\frac{r^2}k\ge\frac r{h_{\max}-1},
\]
so $(\rho_1k_2-\rho_2k_1)^2\le k_1k_2k\cdot r(2C_0+r)/(h_{\min}-r)^2$, that is, $|\rho_1k_2-\rho_2k_1|\le C_2\sqrt{h_{\min}}$. Consequently $k_j\ge k/(2r)$ and $|\pi_1-\pi_2|=\varepsilon|\rho_1k_2-\rho_2k_1|/(k_1k_2)\le C_3h_{\min}^{-3/2}$.

Fix a line $L$, and let $H_L$ be the $a_1(L)$-regular graph of $R_1$-pairs inside $L$. Let $x$ be a unit eigenvector of $H_L$, supported on $L$ and orthogonal to $\chi_L$, with eigenvalue $\vartheta$, and let $x_i=\|E_ix\|^2$. Then $x_0=0$ and $x^{\top}A_1x=\vartheta$, so
\[
\vartheta+\rho_1=(P_{11}+\rho_1)x_1+(P_{21}+\rho_1)x_2,\qquad x_1+x_2=x^{\top}\Pi x=\frac dn-\pi_2+(\pi_1-\pi_2)\vartheta.
\]
Since $|P_{i1}|\le k_1$, substituting the second equation into the first gives
\[
|\vartheta|\le\rho_1+(k_1+\rho_1)\Bigl(\frac dn+\pi_2\Bigr)+(k_1+\rho_1)|\pi_1-\pi_2|\,|\vartheta|.
\]
Here $k_1\le rh_{\max}\le r(h_{\min}+2C_0)$, $d/n\le r/h_{\min}$, $\pi_2\le\rho_2/k_2\le 2r^2/k$, so that $(k_1+\rho_1)\pi_2\le2r^2(k+r)/k$ with $k\ge h_{\min}-1$, and $(k_1+\rho_1)|\pi_1-\pi_2|\le C_3(rh_{\max}+r)h_{\min}^{-3/2}$, which tends to $0$ as $n\to\infty$. Hence, using $(k_1+\rho_1)d/n\le r^2+o(1)$, the second term is at most $3r^2+o(1)$, and moving the last term to the left side gives $|\vartheta|\le C_4:=3r^2+2r$ for $n\ge n_t$. Summing $\vartheta^2$ over the eigenvalues of $H_L$ other than $a_1(L)$ gives $a_1(L)(h_L-a_1(L))\le(h_L-1)C_4^2$. So every line has $a_1(L)\le2C_4^2$ (\emph{type $2$}) or $a_2(L)\le2C_4^2$ (\emph{type $1$}).

Finally, an $R_2$-pair inside a type-$1$ line has at least $h_{\min}-2-4C_4^2$ common $R_1$-neighbors on that line, while an $R_2$-pair inside a type-$2$ line has at most $2C_4^2+(r-1)^2$ common $R_1$-neighbors in total. Since $p^2_{11}$ is a single number, these two kinds of pairs cannot both exist. If there were no type-$2$ line, then $k_2\le2rC_4^2$, contradicting \eqref{eq:large}. Hence type-$1$ lines contain no $R_2$-pairs, and symmetrically type-$2$ lines contain no $R_1$-pairs. So the partial linear space is monochromatic. By the identity $d/n+(h_L-1)\pi_j=1$ for an $R_j$-line $L$, all $R_j$-lines have the same size, so every point lies on the same number $r_j$ of $R_j$-lines.
\end{proof}

From now on we use Claim~\ref{cl:geom}. Let $s_j+1$ be the common size of the $R_j$-lines; as in the last paragraph, $\pi_j=\varepsilon/s_j$ for $j=1,2$. Also $r_1,r_2\ge2$, since a connected union of cliques is complete, and $s_1,s_2\to\infty$ as $n\to\infty$. Let $U_j=\col N_j$, which is the sum of the $E_i$ with $P_{ij}\ne-r_j$. Since $k_j\ne-r_j$ and $P_{3j}=-r_j$, we have $E_0\subseteq U_j$ and $E_3\not\subseteq U_j$; moreover $U_1+U_2=E_3^{\perp}$, and $U_j=E_0$ would give $A_j+r_jI=cJ$, making $R_j$ complete. So, after renaming $E_1$ and $E_2$, one of the following holds:
\begin{itemize}
\item[(a)] $U_1=E_0\oplus E_1$ and $U_2=E_0\oplus E_2$;
\item[(b)] one of $U_1,U_2$ is $E_3^{\perp}$ and the other is $E_0\oplus E_1$;
\item[(c)] $U_1=U_2=E_3^{\perp}$.
\end{itemize}
In each case, if $A_j$ has only the eigenvalues $k_j$, $P_{1j}$ and $-r_j$, then $R_j$ is a connected strongly regular graph with smallest eigenvalue $-r_j$. Its complement is not a disjoint union of cliques, since the scheme is primitive, so $R_j$ is not complete multipartite. By Theorem~\ref{thm:N}, for $n\ge n_t$ it is a Steiner graph or a Latin square graph. Its lines from Claim~\ref{cl:geom} are maximal cliques of $R_j$: a vertex $w$ outside an $R_j$-line $L$ that is $R_j$-adjacent to all of $L$ would lie on $|L|=s_j+1$ distinct $R_j$-lines, one through each point of $L$, while $w$ lies on only $r_j$ lines. Since these cliques have unbounded size, each lies in a point-clique or a net line, which is itself a clique; by maximality the two coincide.

\begin{claim}\label{cl:cases}
None of the cases (a), (b) and (c) occurs for $n\ge n_t$.
\end{claim}

\begin{proof}
\emph{Case (a).} Here $(A_1+r_1I)(A_2+r_2I)$ vanishes on $E_1\oplus E_2\oplus E_3$, so it is a multiple of $J$. No two points are both $R_1$- and $R_2$-collinear, so comparing entries gives $p^1_{12}=r_2(r_1-1)$, $p^2_{12}=r_1(r_2-1)$ and $p^3_{12}=r_1r_2$, the largest possible values. Hence every $R_1$-line meets every $R_2$-line. Indeed, let $L$ be an $R_1$-line and $M$ an $R_2$-line, and choose $u\in L$ and $v\in M$ with $u\ne v$. The entry $\bigl((A_1+r_1I)(A_2+r_2I)\bigr)_{uv}=(N_1N_1^{\top}N_2N_2^{\top})_{uv}$ is the number of pairs $(L',M')$ of an $R_1$-line $L'\ni u$ and an $R_2$-line $M'\ni v$ with $L'\cap M'\ne\emptyset$, since two such lines meet in at most one point. This entry equals $r_1r_2$, the total number of such pairs, so all of them meet; in particular $L\cap M\ne\emptyset$. Moreover, $A_j$ has only the eigenvalues $k_j$, $P_{jj}$ and $-r_j$, so by the remark above both $R_1$ and $R_2$ are Steiner or Latin square graphs, and $n=(s_1+1)(s_2+1)$.

If both are Latin square graphs on $q^2$ points, then the two nets together form an orthogonal array, $G=LS_{r_1+r_2}(q)$ is strongly regular, and $\lambda_1=\lambda_2$. If $R_1$ is the Steiner graph of a $2$-$(m,r_1,1)$ design, then every $R_2$-line meets every point-clique and is an $R_1$-coclique, so it is a parallel class of size $m/r_1$. If moreover $R_2$ is a Latin square graph, then $(m/r_1)^2=m(m-1)/(r_1(r_1-1))$, which forces $m=r_1$. If instead $R_2$ is the Steiner graph of a $2$-$(m',r_2,1)$ design, then the $R_2$-lines are its point-cliques, of size $(m'-1)/(r_2-1)$, and by symmetry the $R_1$-lines are parallel classes of this design, of size $m'/r_2$. Hence
\[
\frac{m'-1}{r_2-1}=\frac m{r_1}\quad\text{and}\quad\frac{m'}{r_2}=\frac{m-1}{r_1-1}.
\]
Eliminating $m'$ gives $(m(r_2-1)+r_1)(r_1-1)=r_1r_2(m-1)$, that is, $m(r_1+r_2-1)=r_1(r_1+r_2-1)$, so again $m=r_1$. Both are impossible, and the remaining case is symmetric.

\emph{Case (b).} Let $U_2=E_0\oplus E_1$ and $U_1=E_3^{\perp}$. Then $A_2$ has the eigenvalues $k_2$, $P_{12}$ and $-r_2$, so $R_2$ is a Steiner graph of a $2$-$(m,K,1)$ design or a Latin square graph $LS_K(q)$, with $K=r_2$, and $E_1$ is its non-trivial positive eigenspace, of dimension $D=m-1$ or $D=K(q-1)$. Write $P_{11}=-\kappa$ and $\beta=P_{21}$.

$A_1$ acts as a scalar on $E_1$. In the Steiner case, $E_1$ consists of the vectors $u\mapsto\sum_{p\in u}f(p)$ over blocks $u$, with $f\perp\mathbf 1$ a function on points. Since $R_1$-adjacent blocks are disjoint, this means that for every block $u$, every point outside $u$ lies in exactly $\kappa$ blocks of $R_1(u)$; so $Kk_1=\kappa(m-K)$ with $\kappa\ge1$. In the Latin square case, applying $A_1$ to $\chi_M-\frac1q\mathbf 1$ for a line $M$ shows that every point off $M$ has exactly $\kappa$ $R_1$-neighbors on $M$, so $k_1=\kappa(q-1)$. In both cases $s_1=k_1/r_1\ge(D-K+1)/(Kr_1)$.

Comparing $A_1^2=k_1I+p^1_{11}A_1+p^2_{11}A_2+p^3_{11}A_3$ on $E_2$ and on $E_3$, where $A_2=-K$ on both, and using $\beta\ne-r_1$, gives $\beta=r_1+p^1_{11}-p^3_{11}$. Since $p^1_{11}\le s_1-1+(r_1-1)^2$, we get $0<\beta+r_1\le s_1+r_1^2$. Taking traces, $m_2(\beta+r_1)=r_1(n-1-D)-k_1+\kappa D$, so $m_2\ge(r_1(n-1-D)-k_1)/(s_1+r_1^2)$. Since $1+D+m_2=\rank N_1\le nr_1/(s_1+1)$,
\[
1+D\le\frac{nr_1(r_1^2-1)}{(s_1+1)(s_1+r_1^2)}+\frac{r_1(1+D)+k_1}{s_1+r_1^2},
\]
which gives $1+D\le2(r_1^3n/s_1^2+r_1)$. As $n\le D^2$ and $s_1\ge(D-K+1)/(Kr_1)$, we get $n/s_1^2\le4K^2r_1^2$ for $D\ge2K$, so $D$, and hence $n$, is bounded.

\emph{Case (c).} From $\Pi^2=\Pi$ on the diagonal,
\[
\frac dn-\Big(\frac dn\Big)^2=k_1\pi_1^2+k_2\pi_2^2+k_3\pi_3^2\ge\varepsilon^2\Big(\frac{r_1}{s_1}+\frac{r_2}{s_2}\Big),
\]
so $d/n\ge\varepsilon(r_1/s_1+r_2/s_2)$. In case (c), $d=\rank N_1=\rank N_2$, so $d/n\le\min\{r_1/(s_1+1),r_2/(s_2+1)\}$. Hence $d/n>2\varepsilon d/n$, that is, $\varepsilon<1/2$ and $d>n/2$. But $d\le nr_1/(s_1+1)<n/2$ for $n\ge n_t$.
\end{proof}

Claim~\ref{cl:cases} shows that the primitive case does not occur for $n\ge n_t$. This completes the proof of Theorem~\ref{thm:main}.
\end{proof}

\begin{proof}[Proof of Corollary~\ref{cor:three}]
It suffices to consider (II). In the proof of Theorem~\ref{thm:main}, case (II) arises in the paragraph after Claim~\ref{cl:W} with $G=S^{\mathcal P}$, where $S=R_2$ has smallest eigenvalue $-a$ and $\tau=a$; by Claim~\ref{cl:W}(i), the smallest eigenvalue of $G$ is $-a-1$, so $a\le t-1$. If $a=1$, then $S$ is a disjoint union of at least two complete graphs (a single complete graph would give $b=1$), so $G$ is a grid, with $p\ne q$ because $G$ has two $\lambda$-values. Let $a\ge2$. Then $S$ is connected, since a disconnected strongly regular graph is a union of complete graphs. By Theorem~\ref{thm:N}, for $n\ge n_t$, $S$ is a complete multipartite graph, a Latin square graph, or a Steiner graph. If $S$ is complete multipartite, its only Hoffman cocliques are its parts, and $G$ is complete. If $S=LS_a(c)$, then a Hoffman coclique has $c$ vertices and meets every line exactly once, so the parts of $\mathcal P$ form a new parallel class; then $G=LS_{a+1}(c)$ is strongly regular, contradicting the assumption that $G$ has two $\lambda$-values. If $S$ is a Steiner graph, its Hoffman cocliques are its parallel classes, and $\mathcal P$ is a resolution of the design.
\end{proof}

\section{Concluding Remarks}\label{sec:remark}

The constant $n_t$ in Theorem~\ref{thm:main} is not explicit. It comes from the theorem of Neumaier, from the bound on $\mu$ of Koolen, Gebremichel, Yang and Yang, and from the requirement that all cliques produced by Lemma~\ref{lem:local} and Theorem~\ref{thm:M} be large compared with these bounds. The density assumption $k<n-2-(t-1)^2/4$ is used only to exclude, in Claim~\ref{cl:types}, the case in which all non-adjacent pairs lie inside the classes of an equivalence relation, so that $G$ contains a spanning complete multipartite subgraph, and to apply the bound on $\mu$.

The assumption that $G$ has two $\lambda$-values cannot be dropped: the constructions of Ge and Koolen~\cite{GK} and of van Dam, Ge and Koolen~\cite{vDGK} give co-edge-regular graphs with four distinct eigenvalues, bounded smallest eigenvalue and at least three $\lambda$-values. For such graphs the relations $R_1$ and $R_2$ no longer form an association scheme together with $\Gbar$, and the spectral arguments of Section~\ref{sec:main}, which rely on the eigenvalues of $A_1$ and $A_2$, are not available.

\section*{Acknowledgments}

The author thanks Hong-Jun Ge and Jack Koolen for helpful discussions during their visit to SUNY Korea in November 2026. The author used Claude (Anthropic) for language refinement, technical editing, computational verification of examples, and symbolic computations. The author is responsible for all mathematical content of this paper.

\end{document}